\documentclass[10pt]{amsart}
\usepackage{amsfonts,amssymb}
\RequirePackage{ifpdf}
\usepackage{amsmath}
\usepackage{color}
\usepackage{pstcol,pst-node}
\usepackage{graphics}
\usepackage{epic}
\usepackage{curves}

\def\appendix#1{
\addtocounter{section}{1} \setcounter{equation}{0}
\renewcommand{\thesection}{\Alph{section}}
\section*{Appendix \thesection\protect\indent\quad
#1}
}

\catcode`\@=11
\def\marginnote#1{}

\newcount\hour
\newcount\minute
\newtoks\amorpm
\hour=\time\divide\hour by60 \minute=\time{\multiply\hour by60
\global\advance\minute by-\hour}
\edef\standardtime{{\ifnum\hour<12 \global\amorpm={am}%
        \else\global\amorpm={pm}\advance\hour by-12 \fi
        \ifnum\hour=0 \hour=12 \fi
        \number\hour:\ifnum\minute<10 0\fi\number\minute\the\amorpm}}
\edef\militarytime{\number\hour:\ifnum\minute<100\fi\number\minute}

\newcommand{\PP}{\mathbb P}

\newcommand{\beq}{\begin{equation}}
\newcommand{\eeq}{\end{equation}}

\newcommand{\ID}{1\!\!1}

\def\be{\begin{equation}}
\def\ee{\end{equation}}
\def\bea{\begin{eqnarray}}
\def\eea{\end{eqnarray}}
\def\<{\langle}
\def\>{\rangle}

\def\nn{\nonumber}

\def\Tr{{\rm Tr}}
\def\one#1{#1^{\raise5pt\hbox{$\scriptstyle\!\!\!\!1$}}\,{}}
\def\two#1{#1^{\raise5pt\hbox{$\scriptstyle\!\!\!\!2$}}\,{}}
\def\onetwo#1{#1^{\raise5pt\hbox{$\scriptstyle\!\!\!\!\!{12}$}}\,{}}

\newtheorem{theorem}{Theorem}[section]

\theoremstyle{definition}
\newtheorem{df}[theorem]{Definition}

\theoremstyle{remark}

\allowdisplaybreaks

\begin{document}
\title[{\it Embedding of the rank 1 DAHA  into $Mat(2,\mathbb T_q)$ and its automorphisms.}]
{Embedding of the rank 1 DAHA  into $Mat(2,\mathbb T_q)$ and its automorphisms.}
\author{Marta Mazzocco}\thanks{\noindent Department of Mathematical Sciences, Loughborough University, Loughborough LE11 3TU, United Kingdom. E-mail: m.mazzocco@lboro.ac.uk.}

\maketitle

\begin{abstract}
In this review paper we show how the Cherednik algebra of type $\check{C_1}C_1$ appears naturally as quantisation of the group algebra of the monodromy group
associated to the sixth Painlev\'e equation. This fact naturally leads to an 
embedding of the Cherednik algebra of type $\check{C_1}C_1$ into $Mat(2,\mathbb T_q)$, i.e. $2\times 2$ matrices with entries in the quantum torus. For $q=1$ this result is equivalent to say that the Cherednik algebra of type $\check{C_1}C_1$ is Azumaya of degree $2$ \cite{O}. 
By quantising the action of the braid group and of the Okamoto transformations on the monodromy group associated to the sixth Painlev\'e equation we study the automorphisms of 
the Cherednik algebra of type $\check{C_1}C_1$ and conjecture the existence of a new automorphism.
Inspired by the confluences of the Painlev\'e equations, we produce similar embeddings for  the confluent Cherednik algebras $\mathcal H_V,\mathcal H_{IV},\mathcal H_{III},\mathcal H_{II}$ and $\mathcal H_{I},$ defined  in \cite{M}. 
\end{abstract}

\phantom{XXX}\hfill {\em For the 60th birthday of Masatoshi Noumi}

\section{Introduction}

The Painlev\'e sixth equation \cite{fuchs, Sch, Gar1} describes the monodromy preserving deformations of a rank $2$ Fuchsian system  with four simple poles $a_1,a_2,a_3$ and $\infty$. The solution of this Fuchsian system is in general a multi-valued analytic vector--function 
in the punctured Riemann sphere $\mathbb P^1\setminus\{a_1,a_2,a_3,\infty\}$ and its multivaluedness is 
described by the so-called monodromy group, i.e. a subgroup of $SL_2(\mathbb C)$ generated by the images $M_1,M_2,M_3$ of the generators of the fundamental group under the anti-homomorphism:
$$
\rho:\pi_1\left( \mathbb P^1\backslash\{a_1,a_2,a_3,\infty\},\lambda_0\right)\to SL_2(\mathbb C).
$$
In this paper, we introduce flat coordinates on a large open sub-set of the set of all possible monodromy groups obtained in this way (see Theorem \ref{th:PVI-mon}). We then obtain a quantisation of the group algebra of the monodromy group by introducing a canonical quantisation for  these flat coordinates. This quantum algebra is isomorphic to 
the Cherednik algebra of type $\check{C_1}C_1$, i.e. the algebra $\mathcal H$ generated by four elements $V_0,V_1,\check{V_0},\check{V_1}$ which satisfy the following relations \cite{Cher,Sa,NS,St}:
\bea
\label{daha1}
(V_0-k_0)(V_0+k_0^{-1})=0\\
\label{daha2}
(V_1-k_1)(V_1+k_1^{-1})=0\\
\label{daha3}
(\check{V_0}-u_0)(\check{V_0}+u_0^{-1})=0\\
\label{daha4}
(\check{V_1}-u_1)(\check{V_1}+u_1^{-1})=0\\
\label{daha5}
\check{V_1}V_1V_0 \check{V_0}=q^{-1/2},
\eea
where  $k_0,k_1,u_0,u_1,q\in\mathbb C^\star$, such that $q^m\neq 1$, $m\in\mathbb Z_{>0}$.

As a consequence we obtain an 
embedding of the Cherednik algebra of type $\check{C_1}C_1$ into $Mat(2,\mathbb T_q)$, i.e. $2\times 2$ matrices with entries in the quantum torus:

\begin{theorem}\label{th-rep-qT}
The map:
\be
\label{rep1}
V_0\to\left(\begin{array}{cc}k_0-k_0^{-1} - \,i e^{-S_3}&-i\, e^{-S_3}\\
k_0^{-1}-k_0+i\, e^{-S_3}+ i\, e^{S_3}&i\, e^{-S_3}\\
\end{array}\right)
\ee
\be
\label{rep2}
V_1\to\left(\begin{array}{cc}k_1-k_1^{-1} - i\, e^{S_2}&k_1-k_1^{-1} - i\, e^{-S_2}-i\, e^{S_2}\\
i\, e^{S_2}&i\, e^{S_2}\\
\end{array}\right)
\ee
\be
\label{rep3}
\check{V_1}\to\left(\begin{array}{cc}0&- \,i e^{S_1}\\
i\, e^{-S_1}&u_1-u_1^{-1}\\
\end{array}\right)
\ee
\be
\label{rep4}
\check{V_0}\to\left(\begin{array}{cc} u_0&0\\
q^{\frac{1}{2}} s& -\frac{1}{u_0}\\
\end{array}\right),
\ee
where $S_1,S_2,S_3$ satisfy the following commutation relations:
\be\label{u0}
[S_{1},S_{2}]=[S_{2},S_3]=[S_3,S_1]=i\pi \hbar,\qquad u_0= -i\, e^{-S_1-S_2-S_3},
\ee
for $q=e^{-i\pi \hbar}$ and
$$
s=\overline k_0 e^{-S_1-S_2}+\overline k_1 e^{-S_1+S_3}+\overline u_1 e^{S_2+S_3}+i\, e^{-S_1-S_2+S_3}+i\, e^{-S_1+S_2+S_3} - u_0,
$$
gives and embedding of $\mathcal H$ into $Mat(2,\mathbb T_q)$. In particular, the images of $V_0,\check{V_0},V_1,\check{V_1}$ in  $GL(2,\mathbb T_q)$ satisfy the relations (\ref{daha1},\dots,\ref{daha4}) and (\ref{daha5}), in the quantum ordering is dictated by the matrix product ordering\footnote{By this we mean that the product $A B$ of two matrices $A, B$ whose entries are in $\mathbb T_q$ is computed by keeping the entries of $A$ on the left matrix of the entries of $B$.}.
\end{theorem}

Note that this result was already proved in \cite{M} (using  a different presentation for $\mathcal H$). The purpose of the current paper is to explain this result in the Painlev\'e context and to draw parallels between the theory of the Painlev\'e equations and the theory of the Cherednik algebra of type $\check{C_1}C_1$. 

In particular, we prove that all the known automorphisms of the Cherednik algebra of type $\check{C_1}C_1$ are a quantisation of the action of the braid group on monodromy matrices proposed in \cite{DM,LT} to describe the analytic continuation of the solutions to the sixth Painlev\'e equation. 

Next we deal with the Okamoto transformations of  the sixth Painlev\'e equation and their action on the monodromy group. By quantisation we conjecture the existence of an automorphism of the Cherednik algebra of type $\check{C_1}C_1$ which acts as follows on the parameters:
$$
(u_1,u_0,k_1,k_0)\to \left(\frac{u_1}{\sqrt{u_1 u_0 k_1 k_0}},\frac{u_0}{\sqrt{u_1 u_0 k_1 k_0}},\frac{k_1}{\sqrt{u_1 u_0 k_1 k_0}},\frac{k_0}{\sqrt{u_1 u_0 k_1 k_0}} \right)
$$
We postpone the computation of the action this automorphism on $V_0,V_1,\check{V_0},\check{V_1}$ to a subsequent publication.

Finally, in \cite{M}, the author introduced  confluent versions of the Cherednik algebra of type $\check{C_1}C_1$ by using a concatenation of Whittaker-type limits similar to those introduced in \cite{Cher1}. In this paper we explain the origin of these confluent  Cherednik algebras from the point of view of the Painlev\'e theory. 
In  \cite{MR} the confluence scheme of the Painlev\'e differential equations was explained in terms of  certain geometric operations giving rise to specific asymptotic limits in the classical coordinates $s_1,s_2,s_3$ and parameters. Here, we quantise these asymptotic limits to obtain asymptotic limits for 
 the quantum coordinates $S_1,S_2,S_3$ and of the parameters $k_0,k_1,u_0,u_1$ (see Fig. 1). By taking these limits  in the matrices (\ref{rep1},\dots\ref{rep4}), we produce new matrices which turn out to provide embeddings for the confluent Cherednik algebras $\mathcal H_V,\mathcal H_{IV},\mathcal H_{III},\mathcal H_{II}$ and $\mathcal H_{I}$ in $Mat(2,\mathbb T_q)$. 

 \begin{figure}[h]
\begin{pspicture}(-6,1)(6,-3.2)
 \rput(-5.5,0){$PVI$}
 \rput(-4.,0){$\begin{array}{c}k_0\to \epsilon,\\ \\
e^{S_3}\to\frac{1}{\epsilon}e^{S_3}\\ \end{array}$}
\psline{->}(-5,0)(-3,0)
\rput(-2.7,0){$PV$}
  \psline{->}(-2.4,0)(-.4,0)
  \rput(0,0){$PIV$}
   \rput(-1.3,0){$\begin{array}{c}k_1\to \epsilon,\\ \\
e^{S_2}\to\frac{1}{\epsilon}e^{S_2}\\ \end{array}$}
  \psline{->}(.4,0)(2.4,0)
  \rput(2.8,0){$PII$}
   \rput(1.3,0){$\begin{array}{c}u_1\to \epsilon,\\ \\
e^{S_1}\to\frac{1}{\epsilon}e^{S_1}\\ \end{array}$}
 \psline{->}(3.2,0)(5.3,0)
  \rput(5.6,0){$PI$}
   \rput(4.3,0){$\begin{array}{c}\\ \\
e^{S_3}\to\frac{1}{\epsilon}e^{S_3}\\ \end{array}$}
 \psline{->}(-2.7,-.5)(-2.7,-2.2)
    \rput(-3.7,-1.5){$\begin{array}{c}
e^{S_3}\to\frac{1}{\epsilon}e^{S_3}\\ \end{array}$}
\rput(-2.7,-2.5){$PIII$}
  \psline{->}(-2.1,-2.5)(.6,-2.5)
  \rput(1.3,-2.5){$PIII^{D_7}$}
   \rput(-1,-2.5){$\begin{array}{c}k_1\to \epsilon, \, u_1\to\epsilon u_1\\ \\
e^{S_1}\to\frac{e^{S_1}}{\epsilon},\, e^{S_2}\to\epsilon e^{S_2}\\ \end{array}$}
  \psline{->}(2,-2.5)(4.5,-2.5)
  \rput(5.2,-2.5){$PIII^{D_8}$}
   \rput(3.4,-2.5){$\begin{array}{c}u_1\to \epsilon,\\ \\
e^{S_1}\to\frac{e^{S_1}}{\epsilon},\, e^{S_2}\to\epsilon e^{S_2}\\ \end{array}$}
 \end{pspicture}
   \caption{The \cite{MR} confluence scheme for the Painlev\'e equations denoted here by $PVI, PV, PIV, PIII, PIII^{D_7}, PIII^{D_8}, PII, PI$ and the corresponding rescaling of the quantum shifted shear coordinates $S_1,S_2,S_3$ such that $\lim_{\hbar \to 0} S_i = s_i+\frac{p_i}{2}$, $i=1,2,3$.   }
   \end{figure}

\section{Flat coordinates for the monodromy group of the sixth Painlev\'e equation}\label{se:PVI}

\subsection{Sixth Painlev\'e equation as isomonodromic deformation equation}
We start by recalling without proof some very well known facts about the sixth  Painlev\'e equation and its relation to the monodromy preserving deformations equations \cite{JMU,MJ1}.

The sixth Painlev\'e equation \cite{fuchs, Sch, Gar1},
\begin{eqnarray}
y_{tt}&=&{1\over2}\left({1\over y}+{1\over y-1}+{1\over y-t}\right) y_t^2 -
\left({1\over t}+{1\over t-1}+{1\over y-t}\right)y_t+\nn\\
&+&{y(y-1)(y-t)\over t^2(t-1)^2}\left[\alpha+\beta {t\over y^2}+
\gamma{t-1\over(y-1)^2}+\delta {t(t-1)\over(y-t)^2}\right],
\end{eqnarray}
describes  the monodromy preserving deformations of a rank $2$ meromorphic connection over 
$\PP^1$ with four simple poles $a_1,a_2,a_3$ and $\infty$ (for example we may choose $a_1=0,\, a_2=t,\, a_3=1$):
\begin{equation}\label{eq:fuchs}
\frac{{\rm d} \Phi}{{\rm d} \lambda} =\sum_{k=1}^3 \frac{A_k(t)}{\lambda-a_k} \Phi,
\end{equation}
where\footnote{For simplicity sake, we are recalling here the main facts about the isomonodromic approach in the case when the parameters $\theta_1,\theta_2,\theta_3$ and $\theta_\infty$ are not integers. This is just a technical restriction, all the results proved in the paper are actually valid also when we lift such restriction.}
\begin{eqnarray}
&&{\rm eigen}(A_i)= \pm\frac{\theta_i}{2}, \quad\hbox{for } i=1,2,3,
\quad A_\infty:=-A_1-A_2-A_3\label{eq:eigen1}\\
&& 
A_\infty =
\left(\begin{array}{cc}\frac{\theta_\infty}{2}&0\\
0&-\frac{\theta_\infty}{2}\\
\end{array}\right),\label{eq:eigen2}
\end{eqnarray}
and the parameters $\theta_i$, $i=1,2,3,\infty$ are related to the PVI parameters by
$$
\alpha=\frac{(\theta_\infty-1)^2}{2},\quad\beta=-\frac{\theta_1^2}{2},\quad \gamma=\frac{\theta_3^2}{2},
\quad \delta=\frac{1-\theta_2^2}{2}.
$$
The precise dependence of the matrices $A_1,A_2,A_3$ on the PVI solution $y(t)$ and its first derivative $y_t(t)$ can be found in \cite{MJ1}.

The solution $\Phi(\lambda)$ of the system (\ref{eq:fuchs}) is a multi-valued analytic function 
in the punctured Riemann sphere $\mathbb P^1\setminus\{a_1,a_2,a_3,\infty\}$ and its multivaluedness is 
described by the so-called monodromy matrices, i.e. the images of the generators of the fundamental group under the anti-homomorphism
$$
\rho:\pi_1\left( \mathbb P^1\backslash\{a_1,a_2,a_3,\infty\},\lambda_0\right)\to SL_2(\mathbb C).
$$
In this paper we fix the base point $\lambda_0$ at infinity and the generators of the fundamental group to be $l_1,l_2,l_3 $, where each $l_i$,  $i=1,2,3$, encircles only  the pole $a_i$ once and  $l_1,l_2,l_3 $ are oriented in such a way that
\be\label{eq:ord-mon}
M_1 M_2 M_3 M_\infty=\ID,
\ee
where $M_\infty=\exp(2\pi i A_\infty)$. 

\subsection{Riemann-Hilbert correspondence and monodromy manifold}
\label{sec:rh}

Let us denote by  ${\mathcal F}(\theta_1,\theta_2,\theta_3,\theta_\infty)$ the moduli space of rank $2$ meromorphic connection over $\mathbb P^1$ with four simple poles $a_1,a_2,a_3,\infty$ of the form (\ref{eq:fuchs}). Let ${\mathcal M}(G_1,G_2,G_3,G_\infty)$  denote the moduli of monodromy representations $\rho$ up to Jordan equivalence, with the local monodromy data of $G_i$'s
prescribed by 
$$ G_i:= \Tr(M_i)= 2\cos(\pi \theta_i),\quad i=1,2,3,\infty. $$
Then the Riemann-Hilbert correspondence 
$$
{\mathcal F}(\theta_1,\theta_2,\theta_3,\theta_\infty)\slash{\Gamma}\to {\mathcal M}(G_1,G_2,G_3,G_\infty)\slash GL_2({\mathbb C}),
$$
where $\Gamma$ is the gauge group \cite{bol}, is defined by associating to each Fuchsian system its monodromy representation class. The representation space 
${\mathcal M}(G_1,G_2,G_3,G_\infty)\slash GL_2({\mathbb C})$ is realised as an affine cubic surface (see \cite{jimbo})
\be\label{eq:friecke}
G_{12}^2+G_{23}^2+G_{31}^2+G_{12}G_{23}G_{31}
-\omega_{3}G_{12}-\omega_{1}G_{23}-\omega_{2}G_{31}+\omega_\infty=0,
\ee
where $G_{12},G_{23},G_{31}$ defined as:
$$G_{ij}= \Tr\left(M_i M_j\right),\qquad i,j = 1,2,3, 
$$
and
$$
 \omega_{ij}:=G_i G_j + G_k G_\infty, \quad k\neq i, j, \qquad
\omega_\infty = G_0^2+G_t^2+G_1^2+G_\infty^2+G_0 G_t G_1 G_\infty-4. $$
This cubic surface is called {\it monodromy manifold} of the sixth Painlev\'e equation and it is equipped with the following Poisson bracket:
\bea\label{eq:pbgij}
&&
\{G_{12},G_{23}\}=G_{12}G_{23}+2 G_{31}-\omega_2 \nn\\
&&
\{G_{23},G_{31}\}=G_{23}G_{31}+2 G_{12}-\omega_3 \\
&&
\{G_{31},G_{12}\}=G_{31}G_{12}+2 G_{23}-\omega_1 \nn
\eea
In \cite{iwa}, Iwasaki proved that the triple 
$(G_{12},G_{23},G_{31})$ satisfying the cubic relation (\ref{eq:friecke}) provides a set of coordinates on a large open subset $\mathcal S\subset{\mathcal M}(G_1,G_2,G_3,G_\infty)$. In the following sub-section, we restrict to such open set.

\subsection{Teichm\"uller theory of the $4$-holed Riemann sphere}

The real slice of moduli space ${\mathcal F}(\theta_1,\theta_2,\theta_3,\theta_\infty)$  of rank $2$ meromorphic connections over $\mathbb P^1$ with four simple poles $a_1,a_2,a_3,\infty$ can be obtained as a quotient of the Teichm\"uller space of the $4$--holed Riemann sphere by the mapping class group. This fact allowed us to use the combinatorial description of the Teichm\"uller space of the $4$--holed Riemann sphere in terms of fat-graphs to produce a good parameterisation of the monodromy manifold of PVI \cite{CM}. In this sub-section we recall the main ingredients of this construction.

We recall that according to Fock~\cite{Fock1}~\cite{Fock2}, 
the fat graph associated to a Riemann
surface  $\Sigma_{g,n}$ of genus $g$ and with $n$ holes is a connected three--valent
graph drawn without self-intersections on $\Sigma_{g,n}$
with a prescribed cyclic ordering
of labelled edges entering each vertex; it must be a maximal graph in
the sense that its complement on the Riemann surface is a set of
disjoint polygons (faces), each polygon containing exactly one hole
(and becoming simply connected after gluing this hole).
In the case of a Riemann sphere $\Sigma_{0,4}$ with $4$ holes, the fat--graph is represented in Fig.2 (the fourth hole is the outside of the graph). 

\begin{figure}
\label{loopinvert}
{\psset{unit=0.5}
\begin{pspicture}(-5,-5)(7,5)
%A3
\newcommand{\PATTERN}[1]{%
\pcline[linewidth=1pt](0.3,0.5)(2,0.5)
\pcline[linewidth=1pt](0.3,-0.5)(2,-0.5)
\psbezier[linewidth=1pt](2,0.5)(3,2)(5,2)(5,0)
\psbezier[linewidth=1pt](2,-0.5)(3,-2)(5,-2)(5,0)
\psbezier[linewidth=1pt](2.8,0)(3.6,0.8)(4,0.7)(4,0)
\psbezier[linewidth=1pt](2.8,0)(3.6,-0.8)(4,-0.7)(4,0)
\rput(1,1.2){\makebox(0,0){$s_{#1}$}}
\rput(4.5,2){\makebox(0,0){$p_{#1}$}}
}
\rput(0,0){\PATTERN{1}}
\rput{120}(0,0){\PATTERN{2}}
\rput{240}(0,0){\PATTERN{3}}
\newcommand{\CURVE}{%
\pcline[linecolor=red, linestyle=dashed, linewidth=1.5pt](0.1,.2)(2.1,.2)
\pcline[linecolor=red, linestyle=dashed, linewidth=1.5pt](0.1,-.2)(2.1,-.2)
\psbezier[linecolor=red, linestyle=dashed, linewidth=1.5pt](2.1,.2)(3,1.7)(4.7,1.8)(4.7,0)
\psbezier[linecolor=red, linestyle=dashed, linewidth=1.5pt]{->}(2.1,-.2)(3,-1.6)(4.7,-1.8)(4.7,0)
}
\rput(0,0){\CURVE}
\rput{120}(0,0){\CURVE}
\psarc[linecolor=red, linestyle=dashed, linewidth=1.5pt](0.1,.2){.4}{210}{270}
\end{pspicture}
}\caption{The fat graph of the $4$ holed Riemann sphere. The dashed geodesic corresponds to $G_{12}$.}
\end{figure}
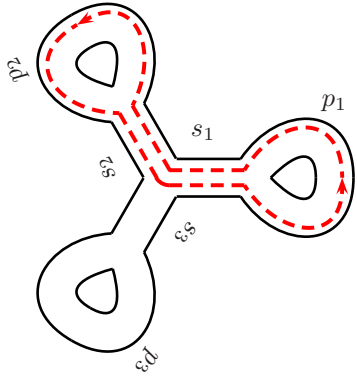

The geodesic length functions, which are traces of hyperbolic elements in the Fuchsian group $\Delta_{g,n}$ such that in Poincar\'e uniformisation:
$$
\Sigma_{g,n}\sim\mathbb H\slash \Delta_{g,n},
$$
are obtained by decomposing each hyperbolic matrix $\gamma\in \Delta_{g,n}$ into a
product of the so--called {\it right, left and edge matrices:}~\cite{Fock1}~\cite{Fock2}
\begin{equation}\nn
R:=\left(\begin{array}{cc}1&1\\-1&0\\
\end{array}\right), \quad
L:=\left(\begin{array}{cc}0&1\\-1&-1\\
\end{array}\right), \quad
E_{s_i}:=\left(\begin{array}{cc}0&-\exp\left({\frac{s_i}{2}}\right)\\
\exp\left(-{\frac{s_i}{2}}\right)&0\end{array}\right).
\label{eq:generators}
\end{equation}
Let us consider the closed geodesics $\gamma_{ij}$ encircling the i-th and j-th holes without self intersections (for example $\gamma_{12}$ is drawn in Fig.1), then their geodesic length functions can be obtained as:
\begin{eqnarray}
\label{eq:shear-PVI-mat}
&&
G_{23}=  -\Tr\left(R E_{s_2} R E_{p_2} R E_{s_2} R E_{s_3} R E_{p_3} R E_{s_3} R\right), \nn\\
&&
G_{31}= - \Tr\left(L E_{s_3} R E_{p_3} R E_{s_3} R E_{s_1} R E_{p_1} R E_{s_1}\right),\\
&&
G_{12}= - \Tr\left( E_{s_1} R E_{p_1} R E_{s_1} R E_{s_2} R E_{p_2} R E_{s_2} L\right),\nn
\end{eqnarray}
which leads to:\footnote{Note that for simplicity we have actually shifted the shear coordinates $s_i\to s_i +\frac{p_i}{2}$, $i=1,2,3$}
\begin{eqnarray}
\label{eq:shear-PVI}
&&
G_{23}=-e^{ s_2+  s_3}-e^{-   s_2- s_3}-e^{-  s_2+  s_3}-G_2e^{ s_3}-G_3 e^{- s_2}\nn\\
&&
G_{31}=-e^{ s_3+ s_1}-e^{- s_3- s_1}-e^{- s_3+ s_1}-G_3e^{ s_1}-G_1e^{- s_3},\\
&&
G_{12}=-e^{ s_1+  s_2}-e^{- s_1-  s_2}-e^{- s_1+  s_2}-G_1e^{ s_2}-G_2 e^{-  s_1}\nn
\end{eqnarray}
where 
$$
G_i=e^{\frac{p_i}{2}}+e^{-\frac{p_i}{2}},\qquad i=1,2,3,
$$
and 
$$
G_\infty=e^{ s_1+ s_2+ s_3}+e^{- s_1- s_2- s_3}.
$$
Due to the classical result that each conjugacy class in the fundamental group $\mathbb P^1\setminus\{a_1,a_2,a_3,\infty\}$ can be represented by a unique closed geodesic, we can make the following identification:
\be\label{eq:identification}
G_{ij}:= \Tr\left(M_i M_j\right),
\ee
and indeed it is a straightforward computation to show that $G_{12},G_{23},G_{31}$ defined as in (\ref{eq:shear-PVI}) indeed lie on the cubic (\ref{eq:friecke}). 

Moreover, the Poisson algebra structure (\ref{eq:pbgij}) 
is induced by the Poisson algebras of geodesic length functions constructed in
\cite{ChF1}  by  postulating the Poisson relations on the level of the
shear coordinates $s_\alpha$ of the Teichm\"uller space. In our case these are:
$$
\{s_1,s_2\}=\{s_2,s_3\}=\{s_3,s_1\}=1,
$$ 
while the perimeters $p_1,p_2,p_3$ are assumed to be Casimirs. 

Since the parameterisation (\ref{eq:shear-PVI}) is analytic in $s_1,s_2,s_3$, we can complexify $s_1,s_2,s_3$ and $G_1,G_2,G_3,G_\infty$ to claim that $s_1,s_2,s_3\in \mathbb C$ provide a system of flat coordinates on the Friecke cubic  (\ref{eq:friecke}).

\subsection{Parameterisation of the monodromy group}

In the case of {\em monodromy matrices}, Korotkin and Samtleben in~\cite{KS}
proposed an $r$-matrix structure of the Fock--Rosly type~\cite{Fock-Rosly}
which did not however satisfy Jacobi relations on monodromy matrices themselves but became
consistent on the level of adjoint invariant elements. Therefore the problem of quantising the monodromy group remained open.

In this section we show that thanks to the identification (\ref{eq:identification}), we can impose
\bea\label{mon-pa}
M_1& = & E_{s_1} R E_{p_1} R E_{s_1},\nn\\
M_2 &=&  - R E_{s_2} R E_{p_2} R E_{s_2} L,\\
M_3 &=&-L E_{s_3} R E_{p_3} R E_{s_3} R,\nn
\eea
so that $s_1,s_2,s_3\in \mathbb C$  provide a system of flat coordinates on the monodromy group, rather than only the monodromy manifold. This will allow as to quantise as we shall see in subsection \ref{suse:q}. 

\begin{theorem}\label{th:PVI-mon}
Given any quadruple of diagonalisable matrices $M_1,M_2,M_3,M_\infty\in SL_2(\mathbb C)$, such that  $M_1 M_2 M_3 M_\infty=\mathbb I$, 
the group $\langle M_1,M_2,M_3\rangle$ is irreducible and none of the matrices $M_1,M_2,M_3,M_\infty$ is a multiple of the identity,
we can find 
 $s_1,s_2,s_3,p_1,p_2,,p_3\in\mathbb C$ such that the following relations hold true (up to global conjugation and cyclic permutation):
\bea\label{mon-par-shear}
M_1&=&\left(\begin{array}{cc}
0&-e^{s_1}\\
e^{-s_1}& - e^{\frac{p_1}{2}}-e^{-\frac{p_1}{2}}\\
\end{array}\right),\nn\\
M_2&=&\left(\begin{array}{cc}
 - e^{\frac{p_2}{2}}-e^{-\frac{p_2}{2}}-e^{s_2}&
  - e^{\frac{p_2}{2}}-e^{-\frac{p_2}{2}}-e^{s_2}-e^{-s_2}\\
  e^{s_2}& e^{s_2}\\
\end{array}\right),\nn\\
M_3&=&\left(\begin{array}{cc}
 - e^{\frac{p_3}{2}}-e^{-\frac{p_3}{2}}-e^{-s_3}&-e^{-s_3}\\
  e^{\frac{p_3}{2}}+e^{-\frac{p_3}{2}}+e^{-s_3}+e^{-s_3}&e^{-s_3}\\
\end{array}\right),\nn\\
M_\infty&=&\left(\begin{array}{cc}
-e^{-s_1-s_2-s_3}&0\\
s_\infty&-e^{s_1+s_2+s_3}
\end{array}\right),\nn\eea
where
\bea
s_\infty&=&\left( e^{\frac{p_3}{2}}+e^{-\frac{p_3}{2}}\right) e^{-s_1-s_2}+\left( e^{\frac{p_2}{2}}+e^{-\frac{p_2}{2}}\right) e^{-s_1+s_3}+
\left( e^{\frac{p_1}{2}}+e^{-\frac{p_1}{2}}\right) e^{s_2+s_3}+\nn\\
&+&e^{-s_1-s_2-s_3}+e^{-s_1-s_2+s_3}+e^{-s_1+s_2+s_3}.\nn
\eea
\end{theorem}

Note that in this parameterisation 
$$
{\rm eigen}(M_j)=-e^{\pm \frac{p_j}{2}},\, j=1,2,3, \hbox{so that }\Tr(M_i)=G_i=e^{\frac{p_i}{2}}+e^{-\frac{p_i}{2}},\qquad i=1,2,3,
$$
and $M_\infty= \left(M_1 M_2 M_3\right)^{-1}$ is not diagonal but has eigenvalues $e^{\pm(s_1+s_2+s_3)}$.

\section{Quantisation}\label{suse:q}

In \cite{CM1}  the proper quantum ordering for a special class of geodesic functions corresponding to geodesics going
around exactly two holes was constructed and it was proved that for each such geodesic, the matrix entries of the corresponding element in the Fuchsian group satisfy a deformed version of  the quantum universal enveloping algebra $U_q(\mathfrak{sl}_2)$ relations. 

In this section we use the same quantum ordering for quantising matrix elements of the monodromy group.

In \cite{CM}, the quantum Painlev\'e cubic can be obtained by introducing the Hermitian operators $S_1,S_2,S_3$ subject to the commutation
inherited from the Poisson bracket of $ s_i$:
$$
[S_i,S_{i+1}]=i\pi \hbar \{ s_i, s_{i+1}\}=i\pi \hbar,\quad i=1,2,3,\ i+3\equiv i,
$$
while the central elements, i.e. perimeters $p_1,p_2,p_3$  and $S_1+ S_2+ S_3$ remain non--deformed, so that the constants $\omega_i^{(d)}$ remain non-deformed  \cite{CM}. 

The Hermitian operators $G_{23}^\hbar,G_{31}^\hbar,G_{12}^\hbar$ corresponding to $G_{23},G_{31},G_{12}$ are introduced as follows: consider the classical expressions for $G_{23},G_{31},G_{12}$ is terms of $ s_1, s_2, s_3$ and $p_1,p_2,p_3$. Write each product of exponential terms as the exponential of the sum of the exponents and replace those exponents by their quantum version, for example the classical $G_{23}$ is
$$
G_{23}=-e^{ s_2+  s_3}-e^{-  s_2- s_3}-e^{-  s_2+  s_3}-G_2e^{ s_3}-G_3 e^{- s_2},
$$
and its quantum version is defined as
$$
G_{23}^\hbar=-e^{S_2+  S_3}-e^{- S_2- S_3}-e^{-  S_2+  S_3}-G_2e^{ S_3}-G_3 e^{- S_2}.
$$
Then $G_{23}^\hbar,G_{31}^\hbar,G_{12}^\hbar$ satisfy the following quantum algebra~\cite{CM}:
\bea
q^{-1/2}G^\hbar_{12}G^\hbar_{23}-q^{1/2}G^\hbar_{23}G^\hbar_{12}&=&(q^{-1}-q)G^\hbar_{13}+(q^{-1/2}-q^{1/2})\omega_{2}\nn\\
q^{-1/2}G^\hbar_{23}G^\hbar_{13}-q^{1/2}G^\hbar_{13}G^\hbar_{23}&=&(q^{-1}-q)G^\hbar_{12}+(q^{-1/2}-q^{1/2})\omega_{3}
\label{q-comm}\\
q^{-1/2}G^\hbar_{13}G^\hbar_{12}-q^{1/2}G^\hbar_{12}G^\hbar_{13}&=&(q^{-1}-q)G^\hbar_{23}+(q^{-1/2}-q^{1/2})\omega_{1}\nn
\eea
and satisfy the following quantum cubic relations:
\bea\label{q-cubics}
{\mathcal C}^\hbar&=&q^{-1/2}G^\hbar_{12}G^\hbar_{23}G^\hbar_{13}-
q^{-1}\bigl(G^\hbar_{12}\bigr)^2
-q \bigl(G^\hbar_{23}\bigr)^2-q^{-1}\bigl(G^\hbar_{13}\bigr)^2\nn\\
&&-q^{-1/2}\omega_{3}G^\hbar_{12}-q^{1/2}\omega_{1}G^\hbar_{23}-q^{-1/2}\omega_{2}G^\hbar_{13},
\eea
where ${\mathcal C}^\hbar$ is a central element is the quantum algebra (\ref{q-comm}).

We now quantise the monodromy matrices in the same way:

\begin{theorem}\label{th:PVI-mon}
The following matrices
\bea\label{mon-q-par}
M_1^\hbar&=&\left(\begin{array}{cc}
0&-e^{S_1}\\
e^{-S_1}& - e^{\frac{p_1}{2}}-e^{-\frac{p_1}{2}}\\
\end{array}\right),\nn\\
M_2^\hbar&=&\left(\begin{array}{cc}
 - e^{\frac{p_2}{2}}-e^{-\frac{p_2}{2}}-e^{S_2}&
  - e^{\frac{p_2}{2}}-e^{-\frac{p_2}{2}}-e^{S_2}-e^{-S_2}\\
  e^{S_2}& e^{S_2}\\
\end{array}\right),\nn\\
M_3^\hbar&=&\left(\begin{array}{cc}
 - e^{\frac{p_3}{2}}-e^{-\frac{p_3}{2}}-e^{-S_3}&-e^{-S_3}\\
  e^{\frac{p_3}{2}}+e^{-\frac{p_3}{2}}+e^{-S_3}+e^{-S_3}&e^{-S_3}\\
\end{array}\right),\nn\\
M_\infty^\hbar&=&\left(\begin{array}{cc}
-e^{-S_1-S_2-S_3}&0\\
s_\infty^\hbar&-e^{S_1+S_2+S_3}
\end{array}\right),\nn\eea
where
\bea
s_\infty^\hbar&=&\left( e^{\frac{p_3}{2}}+e^{-\frac{p_3}{2}}\right) e^{-S_1-S_2}+\left( e^{\frac{p_2}{2}}+e^{-\frac{p_2}{2}}\right) e^{-S_1+S_3}+
\left( e^{\frac{p_1}{2}}+e^{-\frac{p_1}{2}}\right) e^{S_2+S_3}+\nn\\
&+&e^{-S_1-S_2-S_3}+e^{-S_1-S_2+S_3}+e^{-S_1+S_2+S_3},\nn
\eea
are elements of $SL(2,\mathbb T_q)$ and satisfy the following relations:
\bea\label{qmon-algebra}
(M_1^\hbar+ e^{\frac{p_1}{2}}\mathbb I)(M_1^\hbar+ e^{\frac{-p_1}{2}}\mathbb I)=0,\nn\\
(M_2^\hbar+ e^{\frac{p_2}{2}}\mathbb I)(M_2^\hbar+ e^{\frac{-p_2}{2}}\mathbb I)=0,\nn\\
(M_3^\hbar+ e^{\frac{p_2}{2}}\mathbb I)(M_3^\hbar+ e^{\frac{-p_3}{2}}\mathbb I)=0,\nn\\
(M_\infty^\hbar+ e^{S_1+S_2+S_3}\mathbb I)(M_\infty^\hbar+ e^{-S_1-S_2-S_3}\mathbb I)=0,\nn\\
M_\infty^\hbar M_1^\hbar M_2^\hbar M_3^\hbar= q^{-\frac{1}{2}} \mathbb I,
\eea
where $\mathbb I$ is the $2\times 2$ identity matrix.
\end{theorem}
 
This theorem shows that we can interpret the Cherednik algebra as quantisation of the group algebra of the monodromy group of the sixth Painlev\'e equation, in fact the matrices defined by (\ref{rep1}), (\ref{rep2}), (\ref{rep3}), (\ref{rep4}) are simply obtained as $i M_3^\hbar$, $i M_2^\hbar$, $i M_1^\hbar$ and $i M_\infty^\hbar$ respectively so that Theorem \ref{th-rep-qT} can be stated as follows:

\begin{theorem}
The map:
\be
\label{rep1-m}
V_0\to i M_3^\hbar,\quad 
V_1\to i M_2^\hbar,\quad 
\check{V_1}\to i M_1^\hbar,\quad
\check{V_0}\to i M_\infty^\hbar,
\ee
where $M_1^\hbar,M_2^\hbar,M_3^\hbar,M_\infty^\hbar$ are defined as in (\ref{mon-q-par}),
gives and embedding of $\mathcal H$ into $ Mat(2,\mathbb T_q)$. In other words, the matrices $i M_3^\hbar$, $i M_2^\hbar$, $i M_1^\hbar$ and $i M_\infty^\hbar$ 
in  $GL(2,\mathbb T_q)$ satisfy the relations (\ref{daha1},\ref{daha2},\ref{daha3}) and (\ref{daha4}), in which the quantum ordering is dictated by the matrix product ordering and
$$
u_1= - i e^{-\frac{p_1}{2}},\quad k_0 = - i e^{-\frac{p_3}{2}},\quad k_1=-i e^{-\frac{p_2}{2}},\quad
u_0=-i e^{-S_1-S_2-S_3}.
$$
\end{theorem}

\noindent{\it Proof of Theorem \ref{th-rep-qT}.}\/  To prove this theorem, we use the fact that the algebras $\mathcal H$ is the 
 algebra generated by five elements $T^{\pm 1},X^{\pm 1},Y^{\pm 1}$ with the following relations:
\bea
\label{LD0} X W =W X=1,\\
\label{LD00} Y Z =Z Y=1,\\
\label{LD1} X T + a b T^{-1} W +a+b=0,\\
\label{LD2} Z T +\frac{q}{c d} T^{-1} Y + 1+\frac{q}{c d}=0,\\
\label{LD3} (T+ a b)(T+1)=0,\\
\label{LD4} \qquad \quad Y X  =- \frac{q}{a b} T^2 X Y  -q\left( \frac{1}{a}+\frac{1}{b}\right) T Y - \left(1+\frac{c d}{q}\right)T X + (c+d)T,
\eea
where
$$
a=-\frac{u_1}{k_1},\qquad b= u_1 k_1, \qquad c= - \sqrt{q} \frac{k_0}{u_0}, \qquad d=\sqrt{q} u_0 k_0,
$$
and
$$       
X= \sqrt{q} V_0 \check V_0, \qquad Y = k_0 u_1 \check V_1 V_0, \qquad T=u_1 \check V_1,
$$
and viceversa
$$
\check V_1=\frac{1}{u_1} T, \qquad V_0 = \frac{1}{k_0} T^{-1} Y,\qquad \check V_0 = \frac{k_0}{\sqrt{q}}Y^{-1} T X,\qquad V_1=\frac{1}{u_1} T X^{-1}.
$$
We then use Theorem 5.2 from \cite{M} to prove formulae (\ref{rep1},\ref{rep2},\ref{rep3},\ref{rep4}).\hfill{$\square$}
\endproof

 \section{Classical limit of the automorphisms of the Cherednik algebra}

In \cite{DM} the following action of the braid group on monodromy matrices was proposed to describe the analytic continuation of solutions to the sixth Painlev\'e equation:
\bea
&&
\beta_{1}(M_1,M_2,M_3,M_\infty)=( M_1 M_{2}M_1^{-1},M_1,M_3,M_\infty ) ,\nn\\
&&                 
\beta_{2}(M_1,M_2,M_3,M_\infty)=( M_1,M_2 M_{3}M_2^{-1},M_2 ,M_\infty) .        
 \label{braid-M}
\eea
In \cite{LT}, this action was expended  by adding the following involution:
\be\label{braid-r}
r(M_1,M_2,M_3,M_\infty)=(M_3^{-1},M_2^{-1},M_1^{-1},M_\infty^{-1}).
\ee
In this section we prove that this  extended braid group action gives rise to the automorphisms of the Cherednik algebra of type $\check{C_1}C_1$ which were studied in \cite{NS,St}. Here we list them as they appear in \cite{O}:
\begin{gather}
\sigma(\check{V_1},V_1,V_0,\check V_0)=(V_0,V_1,{V}^{-1}_1\check{V}_1{V}_1,V_0\check V_0 V_0^{-1}),\\
\tau(\check{V_1},V_1,V_0,\check V_0)=(\check{V_1},V_1,{V}_0\check{V}_0{V}_0^{-1},{V}_0)\\
\eta(\check{V_1},V_1,V_0,\check V_0)=({V}_1^{-1}\check{V}_1^{-1}{V}_1,{V}_1^{-1},{V}_0^{-1},{V}_0\check{V}_0^{-1}{V}_0^{-1}),\\
\pi(\check{V_1},V_1,V_0,\check V_0)=(\check V_0,\check{V_1},V_1,V_0).
\end{gather}
Indeed by quantising (\ref{braid-M}) and (\ref{braid-r}) we obtain 
\bea
&&
\beta^\hbar_{1}(\check{V_1},V_1,V_0,\check V_0)=( \check{V_1}V_1\check{V_1}^{-1},\check{V_1}, V_0,\check V_0) ,\nn\\
&&                 
\beta^\hbar_{2}(\check{V_1},V_1,V_0,\check V_0)=( \check{V_1},V_1  V_0 V_1^{-1},V_1 ,\check V_0)  ,\nn\\
&&                 
r^\hbar(\check{V_1},V_1,V_0,\check V_0)=( V_0^{-1},V_1^{-1},\check{V_1}^{-1},\check V_0^{-1}).\nn
\eea
It is not hard to check that $\sigma= \beta^\hbar_2\beta^\hbar_1\beta^\hbar_2$, $\tau=\pi^2\beta_1^\hbar \pi^{-2}$ and $\eta=r^\hbar \beta^\hbar_2\beta^\hbar_1\beta^\hbar_2$, so that we can claim that the automorphisms of the Cherednik algebra of type $\check{C_1}C_1$ are indeed the quantisation of the extended modular group action described in \cite{DM,LT}.

The Painlev\'e sixth equation admits also an affine group of bi-rational transformations as described in table 1:
\vskip 2mm

 \begin{center}\begin{tabular}{||c||c|c|c|c||c|c||}
\hline
 & $\theta_1$ &   $\theta_2$  & $\theta_3$  & $\theta_{\infty}$   & $y$ & $t$ \\
   \hline\hline
  $w_1$ & $-\theta_1$ &  $\theta_2$  & $\theta_3$  & $\theta_{\infty}$  &$y$ & $t$ \\
  \hline
  $w_2$ & $\theta_1$ &  $-\theta_2$  & $\theta_3$  & $\theta_{\infty}$   & $y$ & $t$ \\
  \hline
  $w_3$ & $\theta_1$ &  $\theta_2$  & $-\theta_3$  & $\theta_{\infty}$   & $y$ & $t$   \\
  \hline
  $w_{\infty}$ & $\theta_1$ &  $\theta_2$  & $\theta_3$  & $2-\theta_{\infty}$   & $y$ & $t$   \\
   \hline
  $w_{\rho}$ & $\theta_1+\rho$ &  $\theta_2+\rho$  & $\theta_3+\rho$  & $\theta_{\infty}+\rho$ & 
    $y+\frac{\rho}{p}$ & $t$   \\
  \hline\hline
  $r_1$ & $\theta_{\infty}-1$ & $\theta_3$ & $\theta_2$ & $\theta_{1}+1$  &$t/y$ & $t$  \\
  \hline
  $r_2$ & $\theta_3$ &  $\theta_{\infty}-1$ &$\theta_1$ & $\theta_2+1$  & $\frac{t(y-1)}{y-t}$ & $t$   \\
    \hline
  $r_3$ & $\theta_2$  & $\theta_1$ & $\theta_{\infty}-1$& $\theta_3+1$  & $\frac{y-t}{y-1}$ & $t$  \\
   \hline\hline
  $\pi_{13}$ & $\theta_3$ & $\theta_2$ & $\theta_1$ & $\theta_{\infty}$   & $1-y$ & $1-t$   \\
  \hline
  $\pi_{1\infty}$ & $\theta_\infty-1$ & $\theta_2$ & $\theta_3$ & $\theta_1+1$  
  & $1/y$ & $1/t$   \\
    \hline
    $\pi_{12}$ & $\theta_2$ & $\theta_1$ & $\theta_3$ & $\theta_\infty$  
  & $\frac{t-y}{t-1}$ & $\frac{t}{t-1}$  \\
\hline
\label{tb:okam}\end{tabular}\vspace{0.2cm}\\
  \nopagebreak[2] Table 1: Bi-rational transformations for Painlev\'e VI, 
  $\rho=\frac{2-\theta_1-\theta_2-\theta_3-\theta_\infty}{2}$.
  \vspace{0.2cm}\end{center}
  
We note that $w_1,w_2,w_3,w_\infty$ have no effect on the monodromy matrices and therefore on the generators of $\mathcal H$, while $r_1,r_2,r_3$  act as combinations of cyclic permutations and the transformation $r$, while the permutations $\pi_{13},\pi_{1\infty},\pi_{12}$ correspond to a simple re-labelling of the parameters. The only transformation which does not have a simple explanation in terms of monodromy matrices $M_1,M_2,M_3,M_\infty$ is $w_\rho$, which was explained in terms of isomonodromic deformations of a $3\times 3$ linear system with one irregular singularity and one simple pole in \cite{M1}. On the parameters $(u_1,u_0,k_1,k_0)$ this transformation acts as follows: 
$$
(u_1,u_0,k_1,k_0)\to \left(\frac{u_1}{\sqrt{u_1 u_0 k_1 k_0}},\frac{u_0}{\sqrt{u_1 u_0 k_1 k_0}},\frac{k_1}{\sqrt{u_1 u_0 k_1 k_0}},\frac{k_0}{\sqrt{u_1 u_0 k_1 k_0}} \right)
$$
We postpone the computation of the action this automorphism on $V_0,V_1,\check{V_0},\check{V_1}$ to a subsequent publication, this will involve a quantum version of the middle convolution operation discussed in \cite{Fil}.

\section{Embedding of the confluent Cherednik algebras into $Mat(2,\mathbb T_q)$}\label{se:embedding-others}

The confluent limits of the Cherednik algebra of type $\check{C_1}C_1$ were introduced in \cite{M}, in terms of a different presentation which is equivalent to the following (see Theorem 3.2 in \cite{M}):

\begin{df}\label{de:main}
Let $k_1,u_0,u_1,q\in\mathbb C^\star$, such that $q^m\neq 1$, $m\in\mathbb Z_{>0}$. 
The confluent Cherednik algebras $\mathcal H_V,\mathcal H_{IV},\mathcal H_{III},\mathcal H_{II},\mathcal H_{I}$ are the 
 algebras generated by four elements $V_0,V_1,\check{V_0},\check{V_1}$ satisfying the following relations respectively:
 \begin{itemize}
 \item $\mathcal H_V$:
\bea
\label{dahaV1}
V_0^2+ V_0=0,\\
\label{dahaV2}
(V_1-k_1)(V_1+k_1^{-1})=0,\\
\label{dahaV3}
\check{V_0}^2+u_0^{-1}\check{ V_0}=0,\\
\label{dahaV4}
(\check{V_1}-u_1)(\check{V_1}+u_1^{-1})=0,\\
\label{dahaV5}
q^{1/2}\check{V_1}V_1V_0 =\check{V_0}+u_0^{-1},\\
\label{dahaV6}
q^{1/2}\check{V_0}\check{V_1}V_1=V_0+1.
\eea
\item $\mathcal H_{IV}$:
\bea
\label{dahaIV1}
V_0^2+ V_0=0,\\
\label{dahaIV2}
V_1^2+V_1=0,\\
\label{dahaIV3}
\check{V_0}^2+\frac{1}{u_0}\check{ V_0}=0,\\
\label{dahaIV4}
(\check{V_1}-u_1)(\check{V_1}+u_1^{-1})=0,\\
\label{dahaIV5}
q^{1/2}\check{V_1}V_1V_0 =\check{V_0}+u_0^{-1},\\
\label{dahaIV6}
\check{V_0}\check{V_1}V_1=0,\\
\label{dahaIV7}
V_0\check{V_0}=0.
\eea
\item $\mathcal H_{III}$:
\bea
\label{dahaIII1}
V_0^2=0,\\
\label{dahaIII2}
(V_1-k_1)(V_1+k_1^{-1})=0,\\
\label{dahaIII3}
\check{V_0}^2+\frac{1}{\sqrt{q}}\check{ V_0}=0,\\
\label{dahaIII4}
(\check{V_1}-u_1)(\check{V_1}+u_1^{-1})=0,\\
\label{dahaIII5}
q^{1/2}\check{V_1}V_1V_0 =\check{V_0}+\frac{1}{\sqrt{q}}\,\\
\label{dahaIII6}
q^{1/2}\check{V_0}\check{V_1}V_1=V_0.
\eea
\item $\mathcal H_{II}$:
\bea
\label{dahaII1}
V_0^2+ V_0=0,\\
\label{dahaII2}
V_1^2=0,\\
\label{dahaII3}
\check{V_0}^2+\check{ V_0}=0,\\
\label{daha-lim4-3}
(\check{V_1}-u_1)(\check{V_1}+u_1^{-1})=0,\\
\label{dahaII4}
q^{1/2}\check{V_1}V_1V_0 =\check{V_0}+1,\\
\label{dahaII5}
\check{V_0}\check{V_1}V_1=0,\\
\label{dahaII6}
V_0\check{V_0}=0.
\eea
\item $\mathcal H_{I}$:
\bea
\label{dahaI1}
V_0^2+V_0=0,\\
\label{dahaI2}
V_1^2=0,\\
\label{dahaI3}
\check{V_0}^2+\check{ V_0}=0,\\
\label{dahaI4}
\check{V_1}^2+ \check{V_1}=0,\\
\label{dahaI5}
q^{1/2}\check{V_1}V_1V_0 =\check{V_0}+1,\\
\label{dahaI6}
\check{V_0}\check{V_1}=0,\\
\label{dahaI7}
V_0\check{V_0}=0.
\eea
\end{itemize}
\end{df}

All these algebras  $\mathcal H_V,\mathcal H_{IV},\mathcal H_{III},\mathcal H_{II},\mathcal H_{I}$ admit embeddings in $Mat(2,\mathbb T_q)$ (see Theorems 5.2, 5.3, 5.4, 5.5 and 5.5 in \cite{M}). Here we report these embeddings for the generators $\check V_1,V_1,V_0,\check V_0$ in order to clarify the confluence scheme in accordance with Figure 1. Note that in Figure 1, we also have the algebras $\mathcal H_{III}^{D_7}$ and $\mathcal H_{III}^{D_8}$ for which we don't have a Noumi Stokman \cite{NS} representation and for which we can't prove the embedding into $Mat(2,\mathbb T_q)$, so that the geometric explanation behind these algebras remains conjectural.

\begin{theorem}\label{th-rep-qT-PV}
The map:
\be
\label{rep-lim1}
V_0\to \left(\begin{array}{cc}-1&0\\
1+ i\, e^{S_3}&0\\
\end{array}\right)
\ee
\be
\label{rep-lim2}
V_1\to \left(\begin{array}{cc}k_1-k_1^{-1} - i\, e^{S_2}&k_1-k_1^{-1} - i\, e^{-S_2}-i\, e^{S_2}\\
i\, e^{S_2}&i\, e^{S_2}\\
\end{array}\right)
\ee
\be
\label{rep-lim3}
\check{V_1}\to \left(\begin{array}{cc}0&- \,i e^{S_1}\\
i\, e^{-S_1}&u_1-u_1^{-1}\\
\end{array}\right)
\ee
\be
\label{rep-lim4}
\check{V_0}\to \left(\begin{array}{cc} 0&0\\
q^{\frac{1}{2}} s& -\frac{1}{u_0}\\
\end{array}\right),
\ee
where 
$$
s= e^{-S_1-S_2}+(\frac{1}{k_1}-k_1) e^{-S_1+S_3}+(\frac{1}{u_1}-u_1) e^{S_2+S_3}+i\, e^{-S_1-S_2+S_3}+i\, e^{-S_1+S_2+S_3}.
$$
gives and embedding of $\mathcal H_{V}$ into $Mat(2,\mathbb T_q)$. The images of $V_0,\check{V_0},V_1,\check{V_1}$ in  $Mat(2,\mathbb T_q)$ satisfy the relations (\ref{dahaV1}), (\ref{dahaV2}), (\ref{dahaV3}),  (\ref{dahaV4}), (\ref{dahaV5}), (\ref{dahaV6}) in which the quantum ordering is dictated by the matrix product ordering.
\end{theorem}

\begin{theorem}\label{th-rep-qT-PIV}
The map:
\be
\label{rep-piv1}
V_0\to \left(\begin{array}{cc}-1&0\\
1+ i\, e^{S_3}&0\\
\end{array}\right)
\ee
\be
\label{rep-piv2}
V_1\to \left(\begin{array}{cc}-1 - i\, e^{S_2}&-1-i\, e^{S_2}\\
i\, e^{S_2}&i\, e^{S_2}\\
\end{array}\right)
\ee
\be
\label{rep-piv3}
\check{V_1}\to \left(\begin{array}{cc}0&- \,i e^{S_1}\\
i\, e^{-S_1}&u_1-u_1^{-1}\\
\end{array}\right)
\ee
\be
\label{rep-piv4}
\check{V_0}\to \left(\begin{array}{cc} 0&0\\
q^{\frac{1}{2}} s& -\frac{1}{u_0}\\
\end{array}\right),
\ee
where 
$$
s=  e^{-S_1+S_3}+(\frac{1}{u_1}-u_1) e^{S_2+S_3}+i\, e^{-S_1+S_2+S_3}.
$$
gives and embedding of $\mathcal H_{IV}$ into $Mat(2,\mathbb T_q)$. The images of $V_0,\check{V_0},V_1,\check{V_1}$ in  $Mat(2,\mathbb T_q)$ satisfy the relations (\ref{dahaIV1}), (\ref{dahaIV2}), (\ref{dahaIV3}),  (\ref{dahaIV4}), (\ref{dahaIV5}), (\ref{dahaIV7}) in which the quantum ordering is dictated by the matrix product ordering.

\end{theorem}

\begin{theorem}\label{th-rep-qT-PIII}
The map:
\be
\label{rep-piii1}
V_0\to \left(\begin{array}{cc}0&0\\
 i\, e^{S_3}&0\\
\end{array}\right)
\ee
\be
\label{rep-piii2}
V_1\to \left(\begin{array}{cc}k_1-k_1^{-1} - i\, e^{S_2}&k_1-k_1^{-1} - i\, e^{-S_2}-i\, e^{S_2}\\
i\, e^{S_2}&i\, e^{S_2}\\
\end{array}\right)
\ee
\be
\label{rep-piii3}
\check{V_1}\to \left(\begin{array}{cc}0&- \,i e^{S_1}\\
i\, e^{-S_1}&u_1-u_1^{-1}\\
\end{array}\right)
\ee
\be
\label{rep-piii4}
\check{V_0}\to \left(\begin{array}{cc} 0&0\\
q^{\frac{1}{2}} s & -1\\
\end{array}\right),
\ee
where 
$$
s= e^{-S_1-S_2}+(\frac{1}{k_1}-k_1) e^{-S_1+S_3}+(\frac{1}{u_1}-u_1) e^{S_2+S_3}+i\, e^{-S_1-S_2+S_3}+i\, e^{-S_1+S_2+S_3}.
$$
gives and embedding of $\mathcal H_{III}$ into $Mat(2,\mathbb T_q)$. The images of $V_0,\check{V_0},V_1,\check{V_1}$ in  $Mat(2,\mathbb T_q)$ satisfy the relations  (\ref{dahaIII1}), (\ref{dahaIII2}), (\ref{dahaIII3}),  (\ref{dahaIII4}), (\ref{dahaIII5}), (\ref{dahaIII6}), in which the quantum ordering is dictated by the matrix product ordering.
\end{theorem}

\begin{theorem}\label{th-rep-qT-PII}
The map:
\be
\label{rep-pii1}
V_0\to \left(\begin{array}{cc}-1&0\\
1+ i\, e^{S_3}&0\\
\end{array}\right)
\ee
\be
\label{rep-pii2}
V_1\to \left(\begin{array}{cc} - i\, e^{S_2}&-i\, e^{S_2}\\
i\, e^{S_2}&i\, e^{S_2}\\
\end{array}\right)
\ee
\be
\label{rep-pii3}
\check{V_1}\to \left(\begin{array}{cc}0&- \,i e^{S_1}\\
0&-1\\
\end{array}\right)
\ee
\be
\label{rep-pii4}
\check{V_0}\to \left(\begin{array}{cc} 0&0\\
i \sqrt{q} e^{-S_1}e^{S_2}e^{S_3}-\sqrt{q}\left(u_1-\frac{1}{u_1}\right) e^{S_2}e^{S_3} & 1\\
\end{array}\right),
\ee
gives and embedding of $\mathcal H_{II}$ into $Mat(2,\mathbb T_q)$. The images of $V_0,\check{V_0},V_1,\check{V_1}$ in  $Mat(2,\mathbb T_q)$ satisfy the relations  (\ref{dahaII1}), (\ref{dahaII2}), (\ref{dahaII3}),  (\ref{dahaII4}), (\ref{dahaII5}), (\ref{dahaII6}), in which the quantum ordering is dictated by the matrix product ordering.
\end{theorem}

\begin{theorem}\label{th-rep-qT-PI}
The map:
\be
\label{rep-pi1}
V_0\to \left(\begin{array}{cc}-1&0\\
1+ i\, e^{S_3}&0\\
\end{array}\right)
\ee
\be
\label{rep-pi2}
V_1\to \left(\begin{array}{cc} - i\, e^{S_2}&-i\, e^{S_2}\\
i\, e^{S_2}&i\, e^{S_2}\\
\end{array}\right)
\ee
\be
\label{rep-pi3}
\check{V_1}\to \left(\begin{array}{cc}0&- \,i e^{S_1}\\
0&-1\\
\end{array}\right)
\ee
\be
\label{rep-pi4}
\check{V_0}\to \left(\begin{array}{cc} 0&0\\
q^{\frac{1}{2}} e^{S_2}e^{S_3}&1\\
\end{array}\right),
\ee
gives and embedding of $\mathcal H_{I}$ into $Mat(2,\mathbb T_q)$. The images of $V_0,\check{V_0},V_1,\check{V_1}$ in  $Mat(2,\mathbb T_q)$ satisfy the relations  (\ref{dahaI1}), (\ref{dahaI2}), (\ref{dahaI3}),  (\ref{dahaI4}), (\ref{dahaI5}), (\ref{dahaI7}), in which the quantum ordering is dictated by the matrix product ordering.
\end{theorem}
 
\proof The proof of this Theorem is very similar to the proof if Theorem \ref{th-rep-qT-PII}, and is therefore omitted. \endproof

\vskip 2mm \noindent{\bf Acknowledgements.} The author is specially grateful to P. Etingof, T. Koornwinder, M. Noumi  and J. Stokman  for interesting discussions on this subject.


\begin{thebibliography}{99}

\footnotesize\itemsep=0pt

\bibitem{AW}
Askey R., Wilson J.~A., Some basic hypergeometric orthogonal polynomials that generalise Jacobi polynomials, {\it Memoirs of the AMS,}\/ {\bf 319} (1985).


\bibitem{bol} %used
Bolibruch A.A., The 21-st Hilbert problem for linear Fuchsian systems, {\it Developments in mathematics: the Moscow school}\/ Chapman and Hall, London, (1993).


\bibitem{ChF1}
%Chekhov L., Fock V., Talk at St.~Petersburg
%Meeting on Selected Topics in Mathematical Physics (May 26--29, 1997).\\
 Chekhov L., Fock V., {\it A quantum Techm\"uller space}, {\sl Theor. and Math.
Phys.} {\bf120} (1999), 1245--1259, {http://arxiv.org/abs/math.QA/9908165}{math.QA/9908165}.



\bibitem{CM}
Chekhov L., Mazzocco M., Shear coordinate description of the quantised versal unfolding of a $D_4$ singularity, {J. Phys. A: Math. Theor.}\/
{\bf 43}, (2010) 442002, 13 pages.


\bibitem{CM1}
Chekhov L., Mazzocco M., 
Quantum ordering for quantum geodesic functions of orbifold Riemann surfaces, 
in {\it Topology, Geometry, Integrable Systems, and Mathematical Physics: Novikov's Seminar 2012-2014,}\/ 
American Mathematical Society Translations--Series 2, {\bf 234} (2014).


\bibitem{MR}
Chekhov L., Mazzocco M., Rubtsov V., Painlev\'e monodromy manifolds, decorated character varieties and cluster algebras, {\it arXiv:1511.03851}, (2015).


\bibitem{Cher} %%used
Cherednik I., Double affine Hecke algebras, Knizhnik-Zamolodchikov equations and Macdonald's operators, {\it Int. Math. Res. Not.}\/ (1992), no.9:171--180.

\bibitem{Cher1}%used
Cherednik I., Whittaker limits of difference spherical functions, {\it Int. Math. Res. Not.}\/ (2009), no.20:3793--3842.


\bibitem{DM}
Dubrovin B.A., Mazzocco M., Monodromy of certain Painlev\'e-VI transcendents and ref\/lection group,
{\it Invent. Math.} {\bf 141} (2000), 55--147.


\bibitem{EE}
Eshmatov A. and Eshmatov F., Notes on isomorphism between skein algebra and DAHA, {\it private communication}, (2012).
 
\bibitem{EOR}
Etingof P., Oblomkov A. and Rains E., Generalised double affine Hecke algebras of rank 1 and quantised del Pezzo surfaces, 
{\it Adv. Math.,}\/ {\bf 212} (2007), 749--796.


\bibitem{Fil}
Filipuk G.,
On the middle convolution and birational symmetries of the sixth Painlevé equation, 
{\it Kumamoto J. Math.}\/ {\bf 19 }(2006), 15--23. 

\bibitem{Fock1}%used
Fock V.V., Combinatorial description of the moduli space of
projective structures, {http://arxiv.org/abs/hep-th/9312193}{hep-th/9312193}.


\bibitem{Fock2}%used
Fock V.V., Dual Teichm\"uller spaces, {http://arxiv.org/abs/dg-ga/9702018}{dg-ga/9702018}.

\bibitem{Fock-Rosly}
V.~V.~Fock and A.~A.~Rosly,
{\it Moduli space of flat connections as a Poisson manifold},
%Advances in quantum field theory and statistical mechanics: 2nd Italian-Russian collaboration (Como, 1996),
{\it Internat. J. Modern Phys. B} {\bf 11} (1997), no. 26-27, 3195--3206.


\bibitem{fuchs}
Fuchs R., Ueber lineare homogene Differentialgleichungen zweiter Ordnung mit
drei im Endlichen gelegenen wesentlich singul\"aren, {\it Stellen. Math. Ann.,}\/  {\bf 63} (1907)
301--321.

\bibitem{Gar1}
Garnier R., Solution du probleme de Riemann pour les systemes diff\'erentielles
lin\'eaires du second ordre, {\it Ann. Sci. Ecole Norm.. Sup., }\/ {\bf 43} (1926) 239--352.

\bibitem{IIS}
Inaba M., Iwasaki K., Saito M., 
Dynamics of the sixth Painlev\'e equation,
in {\it Th\'eories asymptotiques et \'equations de Painlev\'e,}\/ {\sl S\'emin. Congr.,}\/ {\bf 14} (2006) 103--167.

\bibitem{IT}
Ito T. and Terwillinger P., 
Double affine Hecke algebras of rank $1$ and the $\mathbb Z_3$-symmetric Askey--Wilson relations,
{\it SIGMA}\/ {\bf 6} (2010) 065, 9 pages.


 \bibitem{iwa}
 Iwasaki K., An Area-Preserving Action of the Modular Group on Cubic Surfaces and the Painlev\'e VI Equation, {\it Comm. Math. Phys.}\/ {\bf 242} (2003) 185--219.


\bibitem{jimbo}
Jimbo M., Monodromy Problem and the Boundary Condition for Some Painlev\'e Equations,
{\it Publ. RIMS, Kyoto Univ.,}\/ {\bf 18} (1982) 1137--1161.

\bibitem{JMU} 
Jimbo M., Miwa T. and Ueno K.,
Monodromy preserving deformations of linear ordinary differential
equations with rational coefficients \text{I},
{\it Physica 2D,}\/  {\textbf{2}}, (1981), no. 2, 306--352



\bibitem{KLS} %% USED
Koekoek R., Lesky P., Swarttouw R.~F., Hypergeometric Orthogonal Polynomials and Their $q$-Analogues, {\sl Spinger Monographs in Mathematics,}\/ 
Springer-Verlag, Berlin, (2010).

\bibitem{K1}%% USED
Koornwinder T.~H., The relationship between  Zhedanov's algebra $AW(3)$ and the double affine Hecke algebra in the rank one case, {\it SIGMA}\/ {\bf 3} (2007),  063, 15 pp.



\bibitem{KS}%used
Korotkin D. and Samtleben H.,
Quantization of coset space $\sigma$-models coupled to two-dimensional gravity,
{\it Comm. Math. Phys.}\/ {\bf 190} (1997), no. 2, 411--457.


\bibitem{LT}
Lisovyy O. Tykhyy Yu.,
 Algebraic solutions of the sixth Painlev\'e equation,
 {\it  J. Geom. Phys.}\/ {\bf 85} (2014), 124--163. 


\bibitem{M}%used
Mazzocco M., Confluences of the Painlev\'e equations, Cherednik algebras and q-Askey scheme, {\it arXiv:1307.6140, }\/ (2013).


\bibitem{M1}%used
Mazzocco M., 
Irregular isomonodromic deformations for Garnier systems and Okamoto's canonical transformations, 
{\it J. London Math. Soc. (2),} {\bf 70} (2004), no. 2, 405--419. 

\bibitem{MJ1} 
Jimbo M. and Miwa T.,
Monodromy preserving deformations of linear ordinary differential
equations with rational coefficients \text{II},
{\it Physica 2D,}\/
{\textbf{2}} (1981),
no. 3, 407--448.

\bibitem{NS}%%used
Noumi M., and Stokman J.~V., Askey--Wilson polynomials: an affine Hecke algebraic approach,  {\it Laredo Lectures on Orthogonal Polynomials and Special Functions,}\/ 
{\sl Adv. Theory Spec. Funct. Orthogonal Polynomials, Nova Sci. Publ., Hauppauge, NY}\/ (2004) 111--144.

\bibitem{O}
Oblomkov A., Double affine Hecke algebras of rank $1$ and affine cubic surfaces, {\it IMRN}\/ {\bf 2004}, no.18:877--912.


\bibitem{Sa}%%used
Sahi S. Nonsymmetric Koornwinder polynomials and duality, {\it Ann. of Math. (2)}\/ {\bf 150}, (1999) no1:267--282.



\bibitem{Sch}%used
Schlesinger L.,
Ueber eine Klasse von Differentsial System Beliebliger
\text{Ordnung} mit \text{Festen Kritischer Punkten},
{\it J. fur Math.,}
{\bf 141}, (1912),
96--145.


\bibitem{St}
Stokman J.~V., Difference Fourier transforms for nonreduced root systems, {\it Selecta Math. (N.S.)}\/ {\bf 9} (2003) no3:409--494.

\bibitem{Ter}Terwilliger P., 
The universal Askey-Wilson algebra 
{\it SIGMA}\/ {\bf 7} (2011) 069, 24 pages, arXiv:1104.2813.


\end{thebibliography}
\end{document}